%% file: Subset_Partition7Lists.tex
\documentclass[10pt, a4paper, liststotoc, bibtotoc, headsepline]{article}

\input{useandnews}

\usepackage{geometry}
\usepackage{graphicx}
\usepackage{times}

\usepackage[utf8]{inputenc}

\newcommand{\ustitel}[1]{ \begin{center}{\Large{\textbf{#1}} } \end{center} }

\RequirePackage[ngerman=ngerman-x-latest]{}

\usepackage[automark, headsepline]{scrpage2}  
\ohead[\pagemark]{\pagemark}
\ihead[\rightmark]{\rightmark} 
\ofoot[]{}
\ifoot{}

\usepackage{wrapfig}
\usepackage{graphicx}

\usepackage{amsmath}
\usepackage{amssymb}

\usepackage[Algorithmus, plain]{algorithm}
\usepackage[]{algpseudocode}

\usepackage{natbib}
\setcitestyle{aysep={}}
\makeatletter
\renewcommand{\@biblabel}[1]{(#1)}
\makeatother

\usepackage{setspace} 
\usepackage{colortbl}
\usepackage{tabularx}

\newcolumntype{L}[1]{>{\raggedright\arraybackslash}p{#1}}
\newcolumntype{C}[1]{>{\centering\arraybackslash}p{#1}}
\newcolumntype{R}[1]{>{\raggedleft\arraybackslash}p{#1}}

\makeatletter
\newcommand*\addsubsec{\secdef\@addsubsec\@saddsubsec}
\newcommand*{\@addsubsec}{}
\def\@addsubsec[#1]#2{\subsection*{#2}\addcontentsline{toc}{subsection}{#1}
	\if@twoside\ifx\@mkboth\markboth\markright{#1}\fi\fi
}
\newcommand*{\@saddsubsec}[1]{\subsection*{#1}\@mkboth{}{}}
\makeatother

\usepackage{url}
\newenvironment{abstrakt}{\begingroup \leftskip 2em \rightskip \leftskip \textbf{Abstract}}{\par \endgroup}

\algdef{SE}[DOWHILE]{Do}{DoWhile}{\algorithmicdo}[1]{\algorithmicwhile\ #1}%

\usepackage{pdflscape}

\usepackage{varwidth}

\usepackage{mathtools}

\begin{document}

\thispagestyle{empty}

\ustitel{An efficient Algorithm to partition a Sequence of Integers into Subsets with equal Sums}
	
\thispagestyle{empty}
	
\begin{center}
	{ \large Alexander B\"uchel$\h{1,2}$\!\!, Ulrich Gille\ss en$\h{2}$\!\!, 	Kurt-Ulrich Witt$\h{1,2}$ }
\end{center}

\begin{center} 
	Bonn-Rhein-Sieg University of Applied Sciences\\
	Department of Computer Science \\
	$^1$b-it Applied Science Institute \\
	$^2$Research Group Discrete Mathematics and Optimization (ADIMO) \\
	Sankt Augustin, Germany \\
	alexander.buechel@smail.inf.h-brs.de, ulrich.gillessen@smail.inf.h-brs.de
	kurt-ulrich.witt@h-brs.de \\ 
	
	\medskip
	
	
	February 20$^\text{th}$, 2018
\end{center}

\vspace{0.5cm}

\begin{abstrakt} 
	
	Let $I_n = \set{0,1, \ldots, n}$ and $k, t$ be non-negative integers such that $t \geq n$ and $k \cdot t = {n+1 \choose 2}$. In this paper we present an efficient algorithm which partitions the elements of $I_n$ into $k$ mutually disjoint subsets $T_j$ such that $\cup_{j=1}^n T_j = I_n$ and $\sum_{x \in T_{j}} x = t$ for each $j \in \set{1,2, \ldots, k}$. The algorithm runs in $\Ord\kl{n \cdot \kl{\frac{n}{2k} + \log \frac{n(n+1)}{2k}}}$ time.
	
	\medskip
	
\end{abstrakt}
	 
\section{Introduction}

For $n \in \nat$ let $I_n = \set{0,1, \ldots, n}$, $I_n^+ = I_n - \set{0}$, and $\Delta_n = \sum_{l=0}^{n} l = \frac{n(n+1)}{2}$. We call a collection of $k$ mutually disjoint subsets $T_j \subseteq I_n$ a $(t_1,t_2, \ldots, t_k)$-partition of $I_n$ if 
\begin{itemize}
	\item[(i)] $\sum_{x \in T_{j}} x = t_j$ for $\ejk$,
	
	\item[(ii)] $\sum_{j=1}^{k} t_j = \Delta_n$ and
	
	\item[(iii)] $\cup_{j=1}^n T_j = I_n$. 
\end{itemize}
Given the set $I_n$ and the non-negative integers $t_j$, $\ejk$, with $\sum_{j=1}^{k} t_j = \Delta_n$ the related decision problem $\Pi(n,t_1, \ldots, t_k)$ is to decide, whether there exists a $(t_1, \ldots, t_k)$-partition of $I_n$. \cite{Fu:1992} show, that for $k, l, t \in \nat$ with $0 < l \leq \Delta_{n}$ and $(k-1)t + l + \Delta_{k-2} = \Delta_{n}$ a $(t, t+1, \ldots, t+k-2, l)$-partition of $I_n$ exists. \cite{Chen:2005} prove, that a $(t_{1}, \ldots, t_{k})$-partition of $I_n$ exists, if $\sum_{j=1}^{k} t_{j} = \Delta_{n}$ and $t_{j} \geq t_{j+1}$ for $1 \leq j \leq k-1$ and $t_{k-1} \geq n$ hold. In \cite{Buechel:2016} we present a $0$/$1$-linear program to solve partition problems.

\medskip

In the special case, where $t_j = t = \const$ we call $T_1, \ldots, T_k$ a $(k,t)$-partition of $I_n$. Given $n, k, t \in \nat$ with $t \geq n$ and $\Delta_n = k \cdot t$ the decision problem reduces to the question, whether a $(k,t)$ partition of $I_n$ exists. 
\cite{Straight:1979} show that for all $k, t$ with $\Delta_{n} = k \cdot t$ and $t \geq n$ a partition of $I_n$ exists.  \cite{Ando:1990} withdraw the condition $\Delta_{n} = k \cdot t$ and prove that $I_n$ can be partitioned into $k$ disjoint subsets $T_j$ with $\sum T_j = t$ if and only if $k(2k-1) \leq k \cdot t \leq \Delta_n$. 

\medskip

Where as the cited papers study for which $k$-tuples $(t_1, \ldots, t_k)$-partitions of $I_n$ exist, we are interested in efficient algorithms to determine partitions. In this paper we consider problem instances $\Pi(n,k,t)$ with $t \geq n$ and $\Delta_n = k \cdot t$. In chapter \ref{sectionAlgorithm} we introduce the recursive algorithm $\Pi\mathit{Solve}$ which efficiently determines a partition for each instance $\Pi(n,k,t)$. Before, in chapter \ref{sectionMeander} we present the so called meander algorithm which solves problem instances $\Pi(n,k,t)$, where $n$ is even and $2k$ is a divisor of $n$ or where $n$ is odd and $2k$ divides $n+1$, respectively. The reason is, that $\Pi\mathit{Solve}$ can be stopped, when one of these conditions is reached, and the remaining partition can be determined directly by means of the meander algorithm. In chapter \ref{sectionComplexity} we analyze the run time complexity of $\Pi\mathit{Solve}$. Chapter \ref{sectionConclusion} summarizes the paper and mentiones some ideas to improve $\Pi\mathit{Solve}$. 

\section{Meander Algorithm}\label{sectionMeander}

For $a \in \nn$ and $b \in \nat$ we denote $\teilt{b}{a}$ if $b$ is a divisor of $a$. Given the problem instance $\Pi(n,k,t)$ the meander algorithm applies if $n$ ist even and $\teilt{2k}{n}$ or if $n$ is odd an $\teilt{2k}{n+1}$, respectively. The algorithm distributes the elements of the set $I_n$ into the subsets $T_j$ such that 
\begin{align}
\sum_{x \in T_j} x & =  t, \,\, \ejk \label{glMeander1}
\end{align}

\subsection{Case: $\boldsymbol{n}$ even and $\boldsymbol{\teilt{2k}{n}}$}\label{subsectionnger}

Figure \ref{figMeanderEven} shows the part of the meander algorithm which solves problem instances $\Pi(n,k,t)$ when $n$ even and $2k$ divides $n$.

\medskip

\begin{figure}[h]
	\rule{13cm}{1pt} 
	
	\verb| |\texttt{{\bf meandereven}$(n, k, t)$;} \\ 
	\verb| |\texttt{{\bf input:} $\phantom{i}I_n, k, t$ with $n$ even, $\teilt{2k}{n}$, $t \geq n$, and $\Delta_n = k \cdot t$;} \\
	\verb| |\texttt{{\bf output:} $T_j$ with $\sum_{x \in T_j} x = t, \,\, \ejk$;} \\
	\verb|    |\texttt{$T_1 := \set{0}$, \, $T_j := \emptyset, \,\, 2 \leq j \leq k$;} \\
	\verb|    |\texttt{{\bf for}} $j := 1$ {\bf to} $k$ {\bf do} \\
	\verb|       |\texttt{{\bf for}} $i := 1$ {\bf to} $\frac{n}{2k}$ {\bf do} \\
	\verb|          |\texttt{(1) $T_j : = T_j \cup \set{2ki-(j-1)}$}; \\ 
	\verb|          |\texttt{(2) $T_j := T_j \cup \set{2k(i-1) + j)}$}; \\ 
	\verb|       |\texttt{{\bf endfor}}; \\
	\verb|    |\texttt{{\bf endfor}}; \\
	\verb| |\texttt{{\bf end.}}
	
	\rule{13cm}{.5pt} 
	
	{\bf \caption{\label{figMeanderEven}Meander Algorithm in case $\boldsymbol{n}$ even and $\boldsymbol{\teilt{2k}{n}}$.}}
	
	\rule{13cm}{1pt} 
\end{figure}

\medskip

To prove that the algorithm determines a correct $(k,t)$-partition of $I_n$ we have to show (i) that the set of elements assigned to the subsets $T_j$ in (1) and (2) is equal to $I_n^+$ and (ii) that the resulting partition fulfills condition (\ref{glMeander1}). We will verify (i) in Lemma \ref{lemmaMeandereven1} and (ii) in Lemma \ref{lemmaMeandereven2}. Let
\begin{align}
X_1(n,k) & = \set{2ki - (j-1) \mid 1 \leq i \leq \frac{n}{2k}, \, \ejk } \\
X_2(n,k) & = \set{2k(i-1) + j \mid 1 \leq i \leq \frac{n}{2k}, \, \ejk } 
\end{align}
be the sets of elements of $I_n^+$ which are distributed in assignment (1) or assignment (2), respectively.

{\lemma\label{lemmaMeandereven1} \, Let $\Pi(n,k,t)$ be a problem instance such that $n$ even and $\teilt{2k}{n}$, then $I_n^+ = X_1(n,k) \cup X_2(n,k)$. 

\medskip

{\bf Proof} \, For each $x \in I_n^+$ there exist unambiguously $i, r$ such that 
\begin{align}
x & = 2k(i-1) + r, \, 1 \leq i \leq \frac{n}{2k}, \, 1 \leq r \leq 2k \label{glMeandereven3}
\end{align}
We consider the two  following sets of remainders $r \in I_{2k}^+$: $R_1 = \set{2k - (j-1) \mid \ejk}$ and $R_2 = \set{j \mid \ejk}$. Since $r \in R_1$, if $k+1 \leq r \leq 2k$, it follows $R_1 \cap R_2 = \emptyset$ and $R_1 \cup R_2 = I_{2k}^+$. Thus with respect to (\ref{glMeandereven3}) we get either 
$$
x = 2k(i-1) + 2k - (j-1) = 2ki - (j-1)
$$
or
$$
x = 2k(i-1) + j 
$$
It follows $x \in X_1(n,k) \cup X_2(n,k)$. Hence we have shown $I_n^+ \subseteq  X_1(n,k) \cup X_2(n,k)$.

\medskip

If $x \in X_1(n,k)$, then $k+1 \leq x \leq n$, and if $x \in X_2(n,k)$ then $1 \leq x \leq n-k$. Thus, if $x \in X_1(n,k) \cup X_2(n,k)$, we have $1 \leq x \leq n$, hence $x \in I_n^+$ and thereby $X_1(n,k) \cup X_2(n,k) \subseteq I_n^+$.} \hfill $\B$

{\lemma\label{lemmaMeandereven2} \, Let $\Pi(n,k,t)$ be a problem instance with $n$ even and $\teilt{2k}{n}$, then the output $T_j$, $\ejk$, of $\mathtt{meandereven}(n,k,t)$ fulfills condition (\ref{glMeander1}). 
	
\medskip

{\bf Proof} \, For each $j \in \set{1, \ldots, k}$ we have
\begin{align*}
\sum_{x \in T_j} x & \, = 
\sum_{i=1}^{\frac{n}{2k}} (2ki - (j-1)) + \sum_{i=1}^{\frac{n}{2k}} (2k(i-1) + j) = 2k  \sum_{i=1}^{\frac{n}{2k}} (2i-1) + \dfrac{n}{2k} \\ 
& = 2k \dfrac{n^2}{4k^2} + \dfrac{n}{2k} = \dfrac{n(n+1)}{2k} = t
\end{align*} 
Qed.} \hfill $\B$ 

{\satz\label{satMeandereven} \, $\mathtt{meandereven}(n,k,t)$ 

\medskip

{\bf a)} determines a correct partition of $I_n$ for all problem instances $\Pi(n,k,t)$ with $n$ even and $\teilt{2k}{n}$, and 

\medskip

{\bf b)} runs in $\Ord(n)$ time.

\medskip 

{\bf Proof} \, {\bf a)} follows immediately from Lemma \ref{lemmaMeandereven1} and \ref{lemmaMeandereven2}, and {\bf b)} is obvious.} \hfill $\B$

\subsection{Case: $\boldsymbol{n}$ odd and $\boldsymbol{\teilt{2k}{n \! + \!  1}}$}\label{subsectionnuger}

To solve problem instances $\Pi(n,k,t)$ with $n$ odd and $\teilt{2k}{n+1}$ we adapt slightly the $\mathtt{mean}$-$\mathtt{dereven}$-algorithm (see Fig. \ref{figMeanderOdd}). The correctness of the $\mathtt{meanderodd}$-algorithm can be shown analogously to the proof of the correctnes of the $\mathtt{meandereven}$-algorithm. At this point we define the sets of elements asssigned due to labels (1) and (2) in the $\mathtt{meanderodd}$-algorithm as
\begin{align} 
X_1'(n,k) & = \set{2ki - j \mid  1 \leq i \leq \frac{n+1}{2k}, \, \ejk } \\ 
X_2'(n,k) & = \set{2k(i-1) + (j-1) \mid 1 \leq i \leq \frac{n+1}{2k}, \, \ejk }
\end{align}

\medskip

\begin{figure}[h]
	\rule{13cm}{1pt} 
	
	\verb| |\texttt{{\bf meanderodd}$(n, k, t)$;} \\ 
	\verb| |\texttt{{\bf input:} $\phantom{i}I_n, k, t$ with $n$ odd, $\teilt{2k}{n+1}$, $t \geq n$, and $\Delta_n = k \cdot t$;} \\
		\verb| |\texttt{{\bf output:} $T_j$ with $\sum_{x \in T_j} x = t, \,\, \ejk$;} \\
		\verb|    |\texttt{$T_j := \emptyset, \,\, \ejk$;} \\
			\verb|    |\texttt{{\bf for}} $j := 1$ {\bf to} $k$ {\bf do} \\
	\verb|       |\texttt{{\bf for}} $i := 1$ {\bf to} $\frac{n}{2k}$ {\bf do} \\
	\verb|          |\texttt{(1) $T_j := T_j \cup \set{2ki - j}$}; \\ 
	\verb|          |\texttt{(2) $T_j := T_j \cup \set{2k(i-1) + (j-1)}$}; \\ 
	\verb|       |\texttt{{\bf endfor}}; \\
	\verb|    |\texttt{{\bf endfor}}; \\
	\verb| |\texttt{{\bf end.}}
	
	\rule{13cm}{.5pt} 
	
	{\bf \caption{\label{figMeanderOdd}Meander Algorithm in case $\boldsymbol{n}$ odd and $\boldsymbol{\teilt{2k}{n+1}}$.}}
	
	\rule{13cm}{1pt} 
\end{figure}

\medskip

{\lemma\label{lemmaMeanderodd1} \, Let $\Pi(n,k,t)$ be a problem instance such that $n$ odd and $\teilt{2k}{n+1}$, then $I_n = X_1'(n,k) \cup X_2'(n,k)$. 
	
\medskip
	
{\bf Proof} \, For each $x \in I_n$ there exist unambiguously $i, r$ such that 
\begin{align}
	x & = 2k(i-1) + r, \, 1 \leq i \leq \frac{n+1}{2k}, \, 0 \leq r \leq 2k - 1 \label{glMeanderodd3}
\end{align}
We consider the sets of remainders $r \in I_{2k-1}$: $R_1' = \set{2k - j \mid \ejk}$ and $R_2' = \set{j-1 \mid \ejk} = \set{j \mid \njk-1}$. Since $r \in R_1'$, if $k \leq r \leq 2k-1$, it follows $R_1 \cap R_2 = \emptyset$ and $R_1 \cup R_2 = I_{2k-1}$. Thus with respect to (\ref{glMeanderodd3}) we get  
$$
	x = 2k(i-1) + 2k - j = 2ki - j 
$$
or
$$
	x = 2k(i-1) + (j-1)
$$
respectively. It follows $x \in X_1'(n,k) \cup X_2'(n,k)$. Hence we have shown $I_n \subseteq  X_1'(n,k) \cup X_2'(n,k)$.
	
\medskip
	
If $x \in X_1'(n,k)$, then $k \leq x \leq n$, and if $x \in X_2'(n,k)$ then $0 \leq x \leq n-k$. Thus, if $x \in X_1'(n,k) \cup X_2'(n,k)$, we have $0 \leq x \leq n$, hence $x \in I_n$ and thereby $X_1'(n,k) \cup X_2'(n,k) \subseteq I_n$.} \hfill $\B$

{\lemma\label{lemmaMeanderodd2} \, Let $\Pi(n,k,t)$ be a problem instance with $n$ odd and $\teilt{2k}{n+1}$, then the output $T_j$, $\ejk$, of $\mathtt{meanderodd}(n,k,t)$ fullfills condition (\ref{glMeander1}). 
	
\medskip
	
{\bf Proof} \, For each $j \in \set{1, \ldots, k}$ we have
\begin{align*}
	\sum_{x \in T_j} x & = 
	\sum_{i=1}^{\frac{n+1}{2k}} (2ki - j) + \sum_{i=1}^{\frac{n+1}{2k}} (2k(i-1) + (j-1)) = 2k  \sum_{i=1}^{\frac{n}{2k}} (2i-1) - \dfrac{n+1}{2k} \\ 
	& = 2k \dfrac{(n+1)^2}{4k^2} - \dfrac{n+1}{2k} = \dfrac{n(n+1)}{2k} = t
\end{align*} 
Qed.} \hfill $\B$  

{\satz\label{satMeanderodd} \, $\mathtt{meanderodd}(n,k,t)$ 
	
\medskip

{\bf a)} determines a correct partition of $I_n$ for all problem instances $\Pi(n,k,t)$ with $n$ odd and $\teilt{2k}{n+1}$, and

\medskip

{\bf b)} runs in $\Ord(n)$ time.

\medskip

{\bf Proof} \, {\bf a)} follows from Lemma \ref{lemmaMeanderodd1} and \ref{lemmaMeanderodd2}, and {\bf b)} is obvious.} \hfill $\B$
		
\section{The Algorithm $\boldsymbol{\Pi\mathit{Solve}}$}\label{sectionAlgorithm}

In this section we present the different cases which the $\Pi\mathit{Solve}$-algorithm distinguishes. The input to the algorithm is the set $I_n$ and integers $k, t \in \nat$ with $t \geq n$ and $\Delta_n = k \cdot t$. The output is a partition $T_j$, $\ejk$, which fullfills condition (\ref{glMeander1}). 
	
\subsection{Case: $\boldsymbol{2n > t}$}\label{subsection2n1lt}
		
In this case the algorithm makes a distinction between the cases $t$ even and $t$ odd. 
		
\medskip
		
\subsubsection{Case: $\boldsymbol{t}$ even} 
		
The algorithm starts with filling $\frac{2n-t}{2}$ sets as follows:
\begin{align}
	T_j & = \set{ t - n + (j-1), n - (j-1) }, \, 1 \leq j \leq \dfrac{2n-t}{2} \label{glteven1}
\end{align}
Obviously these sets are disjoint and fullfill condition (\ref{glMeander1}). The union of these sets is the set $\set{t-n, \ldots, \frac{t}{2} - 1, \frac{t}{2} + 1, \ldots, n}$. Thus the elements of the set $I_{t-n-1}$ 
and the element $\frac{t}{2}$ remain, these have to be distributed into the empty $k - \frac{2n-t}{2}$ sets. To do this, each of these sets is split into two subsets:
\begin{align}
	T_j & = T_{j,1} \cup T_{j,2}, \, \frac{2n-t}{2} + 1 \leq j \leq k \label{glteven2}
\end{align}
The total number of these subsets is $2(k-n) + t$. The set $T_{\frac{2n-t}{2}+1,1}$ is filled with the element $\frac{t}{2}$:
\begin{align}
	T_{\frac{2n-t}{2}+1,1} & = \set{ \frac{t}{2} } \label{glteven3}
\end{align}
Thus it remains to distribute the elements of $I_{t-n-1}$ into the $2(k-n) + t - 1$ sets $T_{\frac{2n-t}{2}+1,2}$ and $T_{j,s}$, $\frac{2n-t}{2} + 2 \leq j \leq k$, $s \in \set{1,2}$, i.e. it remains to solve the problem instance $\Pi(n',k',t')$ where
\begin{align}
	n' & = t - n - 1 \label{glns1}\\ 
	k' & = 2(k-n) + t - 1 \label{glks1}\\
	t' & = \frac{t}{2} \label{glts1}
\end{align}
We have to verify that this instance fulfills the input conditions 
\begin{align}
	\Delta_{n'} & = k' \cdot t' \label{glveri03}
\end{align}
and
\begin{align}
	t' & \geq n' \label{glveri031}
\end{align}
Using (\ref{glns1}) -- (\ref{glts1}) we get on one side
\begin{align}
	\Delta_{n'} & = \frac{n'(n'+1)}{2}  = \frac{(t-n-1)(t-n)}{2} = \Delta_n + \frac{t^2 - 2tn - t}{2} \label{glveri13}
\end{align}
and on the other side
\begin{align}
	k' \cdot t' & = (2(k - n) + t - 1) \cdot \dfrac{t}{2} = k \cdot t + \frac{t^2 -2tn - t}{2} \label{glveri23} 
\end{align}
Since for our initial problem $\Pi(n,k,t)$ the condition $\Delta_n = k \cdot t$
holds, the verification of (\ref{glveri03}) follows immediately from (\ref{glveri13}) and (\ref{glveri23}). 

\medskip

From $2n > t$ and $t$ even it follows $2n - 2 \geq t$ and from this we get $\frac{t}{2} \geq t - n - 1$. Using (\ref{glns1}) and (\ref{glts1}) condition (\ref{glveri031}) is verified, too.

\medskip

Thus the algorithm can recursively continue to solve the initial problem by determining a solution for the instance $\Pi(n',k',t')$. 

\subsubsection{Case: $\boldsymbol{t}$ odd} 
		
In this case the algorithm initially fills $\frac{2n-t+1}{2}$ sets as follows:
\begin{align}
	T_j & = \set{ t - n + (j-1), n - (j-1) }, \, 1 \leq j \leq \frac{2n-t+1}{2} \label{gltodd1}
\end{align}
Obviously these sets are disjoint and fullfill condition (\ref{glMeander1}). 
The union of these sets builds the set $\set{t-n, \ldots, n}$. Thus the elements of the set $I_{t-n-1}$ 
remain, these have to be distributed into the empty $k - \frac{2n-t+1}{2}$ sets. Thus the instance $\Pi(n',k',t')$ has to be solved, where
\begin{align}
	n' & = t - n - 1  \label{glns2} \\ 
	k' & = k - \frac{2n-t+1}{2} \label{glks2} \\
	t' & = t \label{glts2} 
\end{align}
To proof that this instance is feasible we have to verify, that the input conditions (\ref{glveri03}) and (\ref{glveri031}) are fulfilled in this case as well.
Using (\ref{glns2}) -- (\ref{glts2}) we get on one side	 
\begin{align}
	\Delta_{n'} & = \frac{n'(n'+1)}{2} = \frac{(t-n-1)(t-n)}{2} = \Delta_n + \frac{t^2 - 2tn - t}{2} \label{glveri11}
\end{align}
and on the other side
\begin{align}
	k' \cdot t' & = \kl{ k - \frac{2n - t + 1}{2} } \cdot t = k \cdot t + \frac{t^2 -  2tn - t}{2} \label{glveri21} 
\end{align}
Since $\Delta_n = k \cdot t$ the verification of (\ref{glveri03}) follows immediately from (\ref{glveri21}) and (\ref{glveri23}).

\medskip

From $2n > t$ it follows $n > t - n - 1$. From this we get by means of the input condition $t \geq n$ and the definitions (\ref{glns2}) and (\ref{glts2}): $t' = t \geq n > t - n - 1 = n'$, i.e. condition (\ref{glveri031}) is fulfilled. 
		
\subsection{Case: $\boldsymbol{2n \leq t}$}\label{subsection2nlesst}
		
In this case each set $T_j$ is split into two disjoint subsets: $T_j = T_{j,1} \cup T_{j,2}$, $\ejk$. The sets $T_{j,1}$ will be filled as follows:
\begin{align}
	T_{j,1} & = \set{ n - 2k + j, n - (j-1) } \label{gl2nt1}
\end{align}
Thereby the elements $n-2k+1, \ldots, n$ are already disributed, and the two elements in  each of these sets add up to
\begin{align}
	n - (i-1) + n - 2k + i & = 2(n-k) + 1
\end{align}
It remains to partition the elements of $I_{n-2k}$ 
into the sets $T_{j,2}$ such that the sum of elements in each $T_{j,2}$ equals $t - (2(n-k)+1)$. Thus it remains to solve the problem instance $\Pi(n',k',t')$ with
\begin{align}
	n' & = n - 2k \label{gl219} \\ 
	k' & = k \label{gl220} \\
	t' & = t - 2(n-k) - 1 \label{glt2k} 
\end{align}
As well as in the former cases we have to assure, that the input conditions (\ref{glveri03}) and (\ref{glveri031}) are fulfilled.
On the one side we have
\begin{align}
	\Delta_{n'} & = \dfrac{(n-2k)(n-2k+1)}{2} = \Delta_n + 2k^2 - k - 2kn \label{glveri12}
\end{align}
and on the other side
\begin{align}
	k' \cdot t' & = k \cdot \kl{ t - 2(n-k) - 1 } = k \cdot t -2kn + 2k^2 - k \label{glveri22}
\end{align}
(\ref{glveri03}) follows immediately from (\ref{glveri12}) and (\ref{glveri22}).

\medskip

From $t \geq 2n$ it follows $n + 1 \geq 4k$. By subtraction we get $t - n - 1 \geq 2n - 4k$ and from this and definitions (\ref{gl219}) and (\ref{glt2k})  $t' = t - 2n + 2k - 1 \geq n - 2k = n'$, i.e. condition (\ref{glveri031}) is verified.

\medskip

The considerations so far lead to the algorithm $\Pi\mathit{Solve}$ shown in figure \ref{figure1}, and we proved that it works correctly in all cases.

\begin{figure}[p]
	\rule{13cm}{1pt}

	\verb| |\texttt{$\boldsymbol{\Pi\mathit{Solve}(n,k,t)}$;} \\
	\verb| |\texttt{{\bf input:} $\phantom{i}I_n, k, t$ with $t \geq n$ and $\Delta_n = k \cdot t$;} \\
	\verb| |\texttt{{\bf output:} $T_j$ with $\sum_{x \in T_j} x = t, \,\, \bigcup_{j=1}^n T_j = I_n, \,\, T_i \cap T_j = \emptyset, \,\, \1 \leq i, j \leq k, \,\, i \not= j$;} \\
	
	\medskip
	
	\hspace*{0.5cm}(1) {\bf case} $\teilt{2k}{n}$ \\ \\ 
	\hspace*{1.05cm} {\bf then} fill $\set{T_j}_{1 \leq j \leq k}$ by $\mathit{meandereven}(n,k,t)$ \\ \\ 
	\hspace*{0.8cm} {\bf case} $\teilt{2k}{n+1}$ \\ \\ 
	\hspace*{1.05cm} {\bf then} fill $\set{T_j}_{1 \leq j \leq k}$ by $\mathit{meanderodd}(n,k,t)$ \\ \\ 
	\hspace*{0.5cm}(2) {\bf case} $t \geq 2n$ \\ \\ 
	\hspace*{1.05cm} {\bf then} {\bf for} $1 \leq j \leq k$ {\bf do} $T_{j,1} = \set{n-2k+j, n-(j-1)}$ {\bf endfor}; \\ \\ 
	\hspace*{2.25cm} fill $\set{T_{j,2}}_{1 \leq j \leq k}$ by $\Pi\mathit{Solve}(n-2k,k,t-2(n-k)-1))$; \\ \\ 
	\hspace*{2.25cm} {\bf for} $1 \leq j \leq k$ {\bf do} $T_{j} = T_{j,1} \cup T_{j,2}$ {\bf endfor} \\ \\ 
	\hspace*{0.5cm}(3) {\bf case} $t < 2n$ {\bf and} $t$ even \\ \\ 
	\hspace*{1.05cm} {\bf then} {\bf for} $1 \leq j \leq \frac{2n-t}{2}$ {\bf do} $T_{j} = \set{t-n+(j-1), n-(j-1)}$ {\bf endfor}; \\ \\ 
	\hspace*{2.25cm} $T_{\frac{2n-t}{2}+1,1} = \set{\frac{t}{2}}$; \\ \\ 
	\hspace*{2.25cm} fill $T_{\frac{2n-t}{2}+1,2}$, $\set{T_{j,1}}_{\frac{2n-t}{2}+2 \leq j \leq k}$ and $\set{T_{j,2}}_{\frac{2n-t}{2}+2 \leq j \leq k}$ \\ \\ 
	\hspace*{2.55cm} by $\Pi\mathit{Solve}(t-n-1,2(k-n)+t-1,\frac{t}{2})$ \\ \\
	\hspace*{2.25cm} {\bf for} $\frac{2n-t}{2} + 1 \leq j \leq k$ {\bf do} $T_{j} = T_{\frac{2n-t}{2} + j,1} \cup T_{\frac{2n-t}{2} + j,2}$ {\bf endfor} \\ \\  
	\hspace*{0.5cm}(4) {\bf case} $t < 2n$ {\bf and} $t$ odd \\ \\ 
	\hspace*{0.9cm} {\bf then for} $1 \leq j \leq \frac{2n-t+1}{2}$ {\bf do} $T_{j} = \set{t-n+(j-1), n-(j-1)}$ {\bf endfor}; \\ \\ 
	\hspace*{2.25cm} fill $\set{T_{j}}_{\frac{2n-t+1}{2}+1 \leq j \leq k}$ by $\Pi\mathit{Solve}(t-n-1,k-\frac{2n-t+1}{2},t)$ \\ \\ 
	{\bf endalgorithm}
	
	\rule{13cm}{.5pt} 
	
	{\bf \caption{\label{figure1}Algorithm $\boldsymbol{\Pi\mathit{Solve}}$.}}
	
	\rule{13cm}{1pt} 
\end{figure}

\medskip

\section{Complexity}\label{sectionComplexity}

In this section we analyse the worst case run time complexity of the $\Pi\mathit{Solve}$-Algorithm. The algorithm consists of four subalgorithms related to the cases we distinguish: (1) $\teilt{2k}{n}$ or $\teilt{2k}{n+1}$,  (2) $t \geq 2n$, (3) $t < 2n$ and $t$ even, (4) $t < 2n$ and $t$ odd. We abbreviate these cases by $m$ (meander), $s$ (smaller), $ge$ (greater even), and $go$ (greater odd),  respectively. Then the run $\Pi\mathit{Solve}(n,k,t)$ can be represented by a sequence $\rho'(n,k,t) \in \set{m, s, ge, go}^{+}$. 

{\bsp \, {\bf a)} Let $n = 1337$. The list of runs for all partitions of $I_{1337}$ is:
\begin{align*}
	\rho'(1337,3,298151) & = m \\ 
	\rho'(1337,7,127779) & = s^{94}\,\mathit{go}\,m \\ 
	\rho'(1337,21,42593) & = s^{30}\,\mathit{go}\,s\,\mathit{ge}\,m \\ 
	\rho'(1337,191,4683) & = ss\,\mathit{go}\,m \\ 
	\rho'(1337,223,4011) & = m \\ 
	\rho'(1337,573,1561) & = \mathit{go}\,m \\ 
	\rho'(1337,669,1337) & = m 
\end{align*}
{\bf b)} Let $n = 9999$, then we have 
\begin{align*}
	\rho'(9999,4444,11250) & = \mathit{ge}\,s^3\,\mathit{ge}^4\,\mathit{go}\,m \\ 
	\rho'(9999,4040,12375) & = \mathit{go}\,s^4\,\mathit{go}\,s^4\,\mathit{go}\,s\,\mathit{ge}\,m \\ 
	\rho'(9999,3960,12625) & = \mathit{go}\,s^{3}\,\mathit{ge}\,\mathit{go}\,s^{8}\,\mathit{go}\,m \\ 
	\rho'(9999,3333,15000) & = \mathit{ge}^3\,\mathit{go}\,m \\ 
	\rho'(9999,12,4166250) & = s^{415}\,\mathit{go}\,s\,\mathit{ge}^2\,m
\end{align*}} \hfill $\B$

\medskip

Let $\alpha$ be a non empty sequence over $\Omega' = \set{m, s, ge, go}$, then $\mathit{first}(\alpha)$ is the first and $\mathit{last}(\alpha)$ the last symbol of $\alpha \in {\Omega'}^+$, and $\mathit{head}(\alpha)$ is the sequence without the last symbol. $\Betr{w}_a$ is the number of occurrences of symbol $a \in \Omega'$ in the sequence $w \in \Omega'^\ast$.

\medskip

\begin{minipage}{13cm}

Obviously we have

{\lemma\label{lemmam} \, Let $\Pi(n,k,t)$ be a problem instance, then $\mathit{last}(\rho'(n,k,t)) = m$ and $m$ is not a member of $\mathit{head}(\rho'(n,k,t))$. } \hfill$\B$

\end{minipage}

\medskip

Thus, we may neglect the last symbol of $\rho'(n,k,t)$ and denote $\rho(n,k,t) = \mathit{head}(\rho'(n,k,t))$. As well we do not need the alphabet $\Omega'$, because $\rho(n,k,t) \in \set{s, \mathit{ge}, \mathit{go}}^{\ast}$. We denote this alphabet by $\Omega$.

\medskip

Next we show, that the last call before the recursion stops with the $m$-case cannot be $s$.

{\lemma\label{lemmasnotlast}  \, Let $\Pi(n,k,t)$ be a problem instance. If $\Betr{\rho(n,k,t)} \geq 1$, then $\mathit{last}(\rho(n,k,t)) \not= s$.
	
\medskip

{\bf Proof} \, We assume $\mathit{last}(\rho(n,k,t)) = s$. Let $\Pi(\nu, \kappa, \tau)$ be the problem instance before the last $s$-call. Then by (\ref{gl219}) and (\ref{gl220}) after the $s$-call we have $\nu' = \nu - 2\kappa$ and $\kappa' = \kappa$. Since the next call is $m$ it has to be $\teilt{2\kappa'}{\nu'}$ or $\teilt{2\kappa'}{\nu' + 1}$, i.e. $\teilt{2\kappa}{\nu - 2\kappa}$ or $\teilt{2\kappa}{\nu - 2k + 1}$. It follows $\teilt{2\kappa}{\nu}$ or $\teilt{2\kappa}{\nu + 1}$. Hence the instance $\Pi(\nu,\kappa, \tau)$ would have been solved by an $m$-call, a contradiction to our assumption $\mathit{last}(\rho(n,k,t)) = s$.} \hfill$\B$

{\kor\label{kornotlast} \,Let $\Pi(n,k,t)$ be a problem instance. If $\Betr{\rho(n,k,t)} \geq 1$, then $\mathit{last}(\rho(n,k,t)) \in \set{\mathit{ge}, \mathit{go}}$. } \hfill $\B$

\subsection{Case: $\boldsymbol{2n > t}$ and $\boldsymbol{t}$ odd}

From $2n > t$ we can conclude $t > 2(t - n - 1)$. Using (\ref{glns2}) and (\ref{glts2}) we get $t' > 2n'$. This leads to

{\lemma\label{lemmagos} \, Let $\Pi(n,k,t)$ be a problem instance with $2n > t$, $t$ odd and $\rho'(n,k,t) = \alpha \mathit{go} \beta$, $\alpha \in \Omega^{\ast}$, $\beta \in \Omega'^+$, then 
	
\medskip

{\bf a)} $\mathit{first}(\beta) = m$, if $\Betr{\beta} = 1$,
	
\medskip

{\bf b)} $\mathit{first}(\beta) = s$, if $\Betr{\beta} \geq 2$. } \hfill$\B$

\medskip

Thus, after the case $\mathit{go}$ the recursion ends by call of the meander algorithm or the recursion continues with the $\mathit{s}$ case either.

{\kor\label{korgos} \, Let $\Pi(n,k,t)$ be a problem instance with $2n > t$ and $t$ odd, then $\Betr{\rho(n,k,t)}_s \geq \Betr{\rho(n,k,t)}_{go}$. } \hfill $\B$

%
%
%
%
%
%
%
%

\subsection{Case: $\boldsymbol{2n > t}$ and $\boldsymbol{t}$ even}

From (\ref{glts1}) it follows immediately
%
\begin{align}
	\Betr{\rho(n,k,t)}_{ge} \leq \log t = \log \frac{n(n+1)}{2k} \label{gllogt}
\end{align}

\subsection{Case: $\boldsymbol{2n \leq t}$}\label{subsection2nlesstc}

In this case if the algorithm performs the instance $\Pi(n,k,t)$, then the next instance to solve may be $\Pi(n',k,t')$ with $n' = n - 2k$ and $t' = t - 2(n-k)-1$ (cf. subsection \ref{subsection2nlesst}, equations (\ref{gl219}) and (\ref{glt2k}), respectively). By $n^{(\ell)}$ and $t^{(\ell)}$ we denote the value of $n$ and $t$ in the $\ell^{\text{th}}$ recursion call in the case $2 n^{(\ell)} \leq t^{(\ell)}$. Thus we have $n^{(0)} = n$, $n^{(2)} = n' = n - 2k$ and $t^{(0)} = t$, $t^{(1)} = t' = t - 2(n-k) - 1$, for example. By induction we get
\begin{align}
	n^{(\ell)} & = n - 2k \cdot \ell \label{gl2kell} \\ 
	t^{(\ell)} & = t - 2n \cdot \ell + 2k \cdot \ell^2 - \ell \nonumber \\
			   & = t - (2(n - k \cdot \ell) + 1) \cdot \ell \label{gl2tell}
\end{align}
Now we determine the order of the maximum value of $\ell$ guaranteeing the condition $2 n^{(\ell)} \leq t^{(\ell)}$. Using (\ref{gl2kell}) and (\ref{gl2tell}) we get
\begin{align*}
	0 & \leq t^{(\ell)} - 2 n^{(\ell)} \\
	  & = t - (2(n - k \cdot \ell) + 1) \cdot \ell - 2(n - 2k \cdot \ell)
\end{align*}
To determine $\ell$ we solve the quadratic equation
\begin{align*}
	0 & = \ell^2 + \dfrac{4k - 2n - 1}{2k} \cdot \ell + \dfrac{t - 2n}{2k}
\end{align*}
which has the solutions
\begin{align*}
	\ell_{1,2} & = - \dfrac{4k - 2n - 1}{4k} \pm \sqrt{\kl{\dfrac{4k - 2n - 1}{4k}}^2 - \dfrac{t - 2n}{2k}} \\ \\ 
			   & = - \dfrac{4k - 2n - 1}{4k} \pm \dfrac{4k - 1}{4k}
\end{align*}
i.e.
\begin{align*}
	\ell_1 & = \dfrac{n}{2k}, \,\,\, \ell_2 = \dfrac{n+1}{2k} - 2
\end{align*}
Finally we get
\begin{align}
	\ell \leq \dfrac{n}{2k}
\end{align}
Thus, we have just proven

{\lemma\label{lemma2kkt} \,Let $\Pi(n,k,t)$ be a problem instance. If $\rho(n,k,t) = s^\ell\,x$ with $x \in \set{\mathit{ge}, \mathit{go}}$, then $\ell \leq \dfrac{n}{2k}$. }\hfill$\B$

\medskip

\begin{minipage}{13cm}

Corollary \ref{korgos}, inequality (\ref{gllogt}) and Lemma \ref{lemma2kkt} lead to

{\satz\label{satzconclusion} \, Let $\Pi(n,k,t)$ be a problem instance.
	
\medskip

{\bf a)} \, Then the recursion depth of $\Pi\mathit{Solve}(n,k,t)$ is $\Ord\kl{\frac{n}{2k} + \log \frac{n(n+1)}{2k}}$.

\medskip

{\bf b)} \, Since the complexity of operations the algorithm performs in each recursion call (assig\-ning elements of $I_n$ to some set $T_j$, arithmetic comparisons and oprerations) is $\Ord(n)$ it follows that the worst case run time complexity of $\Pi\mathit{Solve}$ is 
$$
	\Ord\kl{n \cdot \kl{\frac{n}{2k} + \log \frac{n(n+1)}{2k}}}
$$. \hfill$\B$ }

\end{minipage}

\section{Conclusion}\label{sectionConclusion}
		
In section \ref{sectionAlgorithm} we present the recursive algorithm $\Pi\mathit{Solve}$ which solves following subset partition problems: Given $n, k, t \in \nat$ and $I_n$ with $t \geq n$ and $\Delta_n = k \cdot t$, then the algorithm partitions the integers from $0$ to $n$ into $k$ mutually disjoint sets such that the elements in each set add up to $t$. The recursion can be stopped, if $n$ ist even and ${2k}$ is a divisor ${n}$ or if $n$ is odd and ${2k}$ is a divisor of ${n+1}$, respectively, because in these cases the meander algorithms presented in section \ref{sectionMeander} can be applied, which directly determines a partition.

\medskip

We prove that the algorithm works correctly and runs in 
$$
	\Ord\kl{n \cdot \kl{\frac{n}{2k} + \log \frac{n(n+1)}{2k}}}
$$ 
time for each problem instance $\Pi(n,k,t)$.

\medskip

In \cite{Jagadish:2015} an approximation algorithm for the cutting sticks-problem is presented. Because the cutting sticks-problem can be transformed into an equivalent partitioning problem our algorithms can be applied to the corresponding cutting sticks-problems.

\medskip

Further research may investigate whether ideas from the previous chapters and cited papers can be used to improve the efficiency of the $\Pi\mathit{Solve}$-algorithm. In \cite{Buechel:2016}, \cite{Buechel:2017a} and \cite{Buechel:2017b} we present efficient solutions for problem instances $\Pi(n,k,t)$, where $n = q \cdot k$, $q, k$ odd; $n = m^2 - 1$, $m \geq 3$; $n = p - 1$, $\pinprim$; $n = p$, $\pinprim$; $n = 2p$; and  $\pinprim$, where $\prim$ is the set of prime numbers. Thus we may augment the $\Pi\mathit{Solve}$-algorithm by related conditions to stop further recursion calls. 

\medskip

{\bf Acknowledgement:} We would like to thank Arkadiusz Zarychta who created a tool by means of which we are able to test the algorithm and to analyse experimentally its performance.

	\end{document}

%% file: useandnews.tex
\usepackage{epic}   
\usepackage{eepic}  
\usepackage{epsf}
\usepackage{latexsym}
\usepackage{pifont}
\usepackage{amssymb}
\usepackage{amsmath}
\usepackage{makeidx}
\usepackage{eurosym}
\usepackage{color}
\usepackage{colonequals}
\usepackage{times}
\usepackage{a4}
\usepackage{framed}
\usepackage{lscape}
\usepackage{pict2e}

\usepackage{epsfig}
\usepackage{epstopdf}

\usepackage{hyperref}

\definecolor{logo}{rgb}{0.3,0.5,1.0}

\usepackage{epsfig}
\usepackage{epstopdf}

\numberwithin{equation}{section}

\newcommand{\nat}{\mathbb{N}}
\newcommand{\nn}{\mathbb{N}_{_{0}}}

\newcommand{\prim}{\mathbb{P}}
\newcommand{\pinprim}{p \in \mathbb{P}}

\newcommand{\ejk}{1 \leq j \leq k}
\newcommand{\njk}{0 \leq j \leq k}

\newcommand*\Betr[1]{\left|#1\right|} 

\newcommand{\set}[1]{\left\{ #1 \right\}}   
\newcommand{\kl}[1]{\left(\, #1 \,\right)}   

\newcommand*\vr[1]{\mathbf{#1}} 

\newcommand{\1}{\vr{1}} 

\newcommand*\h[1]{^{#1}}

\newcommand{\Ord}{{\cal O}}

\newcommand{\const}{\sf const}

\newcommand*\teilt[2]{#1 \! \left| #2 \right.}

\newcommand{\B}{{\Box}}

\newtheorem{bsp}{Example}[section]

\newtheorem{satz}{Theorem}[section]

\newtheorem{kor}{Corollary}[section]

\newtheorem{lemma}{Lemma}[section]

\newcounter{mt}[subsubsection]

\newcounter{exercise}[section]

\newcommand{\mytheexercise}{\arabic{section}.\arabic{exercise}}

{\endMakeFrame}
